\documentclass{amsproc}
\usepackage{amsfonts}

\theoremstyle{plain}

\newtheorem{corollary}{Corollary}

\newtheorem{lemma}{Lemma}

\newtheorem{proposition}{Proposition}
\newtheorem{remark}{Remark}

\newtheorem{theorem}{Theorem}
\numberwithin{equation}{section}

\begin{document}
\title[Fixed Point Indices and Invariant Periodic Sets]{Fixed Point Indices
and Invariant Periodic Sets of Holomorphic Systems}
\author{Guang Yuan Zhang}
\address{Department of Mathematical Sciences, Tsinghua University, Beijing
100084, The People's Republic of China}
\email{gyzhang@math.tsinghua.edu.cn\\
gyzhang@mail.tsinghua.edu.cn}
\date{}
\thanks{The author is supported by Chinese NSFC 10271063 and 10571009 }

\begin{abstract}
This note presents a method to study center families of periodic orbits of
complex holomorphic differential equations near singularities, based on some
iteration properties of fixed point indices. As an application of this
method, we will prove Needham's theorem in a more general version.
\end{abstract}

\maketitle

\section{Fixed point indices of holomorphic mappings}

Let $\mathbb{C}^{n}$ be the complex vector space of dimension $n$, let $U$
be an open set in $\mathbb{C}^{n}$ and let $f:$ $U\rightarrow \mathbb{C}^{n}$
be a holomorphic mapping. If $p\in U$ is an isolated zero of $f,$ say, there
exists a neighborhood $V$ with $p\in V\subset \overline{V}\subset U$ such
that $p$ is the unique solution of the equation $f(x)=0$ in $\overline{V}.$
Then we can define the zero index of $f$ at $p$ by

\begin{equation*}
\pi _{f}(p)=\#\{x\in V;f(x)=q\},
\end{equation*}%
where $q$ is a regular value of $f$ such that $|q|$ is small enough and $\#$
denotes the cardinality. $\pi _{f}(p)$ is well defined (see [\ref{Mi2}] or [%
\ref{Zh2}] for the detail).

If $f:$ $U$ $\rightarrow $ $\mathbb{C}^{n}$ is a holomorphic mapping and $p\
$is an isolated fixed point of $f,$ then there is a ball $B$ in $U$ centered
at $p$ so that $p$ is the unique fixed point of $f$ in $\overline{B},$ in
other words, $p$ is the unique zero of the mapping
\begin{equation*}
f-I:\overline{B}\rightarrow \mathbb{C}^{n},
\end{equation*}%
which puts each $x\in \overline{B}$ into $f(x)-x,$ and then the \textit{%
fixed point index} of $f$ at $p$ is well defined by%
\begin{equation*}
\mu _{f}(p)=\pi _{f-I}(p).
\end{equation*}

Fixed point indices of holomorphic mappings have the geometric properties
stated in the following lemmas (see [\ref{Zh1}] and [\ref{Zh2}]). We will
denote by $\Delta ^{n}$ a ball in $\mathbb{C}^{n}$ centered at the origin.

\begin{lemma}
\label{lem4}Let $f:\Delta ^{n}\rightarrow \mathbb{C}^{n}$ be a holomorphic
mapping such that $0\in \Delta ^{n}$ is an isolated fixed point of $f.$ Then
$\mu _{f}(0)\geq 1,$ and the equality holds if and only if the Jacobian
matrix $f^{\prime }(0)$ of $f$ at $0$ has no eigenvalue equal to $1.$
\end{lemma}

\begin{lemma}
\label{Lem3}Let $U$ be an open and bounded set in $\mathbb{C}^{n}$ and let $%
f:\overline{U}\rightarrow \mathbb{C}^{n}$ be a holomorphic mapping such that
$f$ has no fixed point on the boundary $\partial U.$ Then $f$ has only
finitely many fixed points in $U$.
\end{lemma}

\begin{lemma}
\label{lem5}Let $f:\overline{\Delta ^{n}}\rightarrow \mathbb{C}^{n}$ be a
holomorphic mapping such that $0$ is the unique fixed point of $f$ in $%
\overline{\Delta ^{n}}.$ Then for any holomorphic mapping $g:\overline{%
\Delta ^{n}}\rightarrow \mathbb{C}^{n}$ that is sufficiently close to $f$ on
the boundary $\partial \Delta ^{n},$ $g$ has only finitely many fixed points
in $\Delta ^{n}$ and%
\begin{equation*}
\mu _{f}(0)=\sum_{\substack{ g(x)=x  \\ x\in \Delta ^{n}}}\mu _{g}(x).
\end{equation*}
\end{lemma}

This lemma gives us an intuitive interpretation of the fixed point index:
with the same condition of this lemma, for any holomorphic mapping $g:%
\overline{\Delta ^{n}}\rightarrow \mathbb{C}^{n}$ such that $g$ is
sufficiently close to $f$ on the boundary $\partial \Delta ^{n}$ and all
fixed points of $g$ in $\Delta ^{n}$ are simple, $g$ has exactly $\mu
_{f}(0) $ distinct fixed points in $\Delta ^{n}.$ A fixed point $x$ of $g$
is called \emph{simple} if $\det (g^{\prime }(x)-I)\neq 0,$ say, the
Jacobian matrix $g^{\prime }(x)$ has no eigenvalue $1,$ where $I$ is the
unit matrix.

Let $f:\Delta ^{n}\rightarrow \mathbb{C}^{n}$ be a holomorphic mapping such
that $0$ is a fixed point of $f.$ Then for any positive integer $k,$ the $k$%
-th iteration $f^{k}$ of $f$ is well defined in some neighborhood $V_{k}$ of
$0$, where $f^{1}=f,$ $f^{2}=f\circ f,$ $f^{3}=f\circ f\circ f,$ and so on.
The following result can be found in a more general and more precise version
in [\ref{Zh1}] when $n=2,$ and when $n>2$ the proof is similar and will be
given in the appendix section.

\begin{proposition}
\label{Lem2}Let $m>1$ be a prime number and let $f:\Delta ^{n}\rightarrow
\mathbb{C}^{n}$ be a holomorphic mapping such that $0\in \Delta ^{n}$ is an
isolated fixed point of both $f$ and $f^{m}.$ If $f^{\prime }(0)$ has an
eigenvalue that is a primitive $m$-th root of $1$, then%
\begin{equation}
\mu _{f^{m}}(0)>m.  \label{3}
\end{equation}
\end{proposition}

A complex number $\lambda $ is called a primitive $m$-th root of $1,$ if $%
\lambda ^{m}=1$ but $\lambda ^{k}\neq 1$ for $k=1,\dots ,m-1.$ The following
lemma is due to Shub, M. and Sullivan, D. [\ref{SS}].

\begin{lemma}
\label{lem6}Let $m>1$ be a positive integer and let $\Theta :\Delta
^{n}\rightarrow \mathbb{C}^{n}$ be a holomorphic mapping with an isolated
fixed point at the origin $0\in \Delta ^{n}.$ Assume that for each
eigenvalue $\lambda $ of $\Theta ^{\prime }(0),$ either $\lambda =1$ or $%
\lambda ^{m}\neq 1$. Then $0$ is still an isolated fixed point of $\Theta
^{m}$.
\end{lemma}

In the rest of this section we will consider a holomorphic system%
\begin{equation}
\dot{x}=\frac{dx}{dt}=F(x),x\in \Delta ^{n},  \label{4}
\end{equation}%
where $F:\Delta ^{n}\rightarrow \mathbb{C}^{n}$ is a holomorphic mapping
such that $0$ is an isolated zero of $F$, say, $0$ is an isolated
singularity of the system. The following result is well known.

\begin{lemma}
\label{C2} For any $\tau >0,$ there is a ball $B\subset \Delta ^{n}$
centered at the origin $0$ such that the local flow $\phi (t,x)$ of (\ref{4}%
) is real analytic on $[0,\tau ]\times B,$ holomorphic with respect to $x,$ $%
\phi ([0,\tau ]\times B)\subset \Delta ^{n}$ and, putting $\Phi _{\tau
}(x)=\phi (\tau ,x)$, the Jacobian matrices $\Phi _{\tau }^{\prime }(0)$ and
$F^{\prime }(0)$ of $\Phi _{\tau }$ and $F$ at $0,$ respectively, are
related by
\begin{equation}
\Phi _{\tau }^{\prime }(0)=e^{\tau F^{\prime }(0)}.  \label{k1}
\end{equation}
\end{lemma}

Now, we can prove the following result as a consequence of Proposition \ref%
{Lem2}.

\begin{corollary}
\label{C1}Consider the holomorphic system (\ref{4}) and assume the Jacobian
matrix $F^{\prime }(0)$ has an eigenvalue $\omega i,$ with $\omega \neq 0$
being real. Let $\Phi (x)=\phi (\frac{2\pi }{|\omega |},x)$ be the time $%
\frac{2\pi }{|\omega |}$ mapping of the local flow $\phi (t,x)$ of the
system. Then $0$ is an accumulated fixed point of $\Phi .$
\end{corollary}

\begin{proof}
By Lemma \ref{C2}, $\Phi $ is well defined and holomorphic in a neighborhood
of $0.$ To prove the conclusion by contradiction, we assume that $0$ is an
isolated fixed point of $\Phi .$ Then $\mu =\mu _{\Phi }(0)$ is well defined.

Let $m$ be an integer with $m>\mu _{\Phi }(0).$ Then $0$ is an isolated
fixed point of both $\Psi (x)=\phi (\frac{2\pi }{m|\omega |},x)$ and the $m$%
-th iteration $\Psi ^{m}$ of $\Psi ,$ for
\begin{equation*}
\Psi ^{m}(x)=\phi (\frac{m2\pi }{m|\omega |},x)=\Phi (x)
\end{equation*}%
in a neighborhood of $0.$ The previous equality also implies that%
\begin{equation}
\mu _{\Psi ^{m}}(0)=\mu _{\Phi }(0)=\mu .  \label{k2}
\end{equation}

On the other hand, it follows from (\ref{k1}) that $\Psi ^{\prime }(0)$ has
an eigenvalue equal to $e^{\frac{2\pi \omega }{m|\omega |}i},$ which is a
primitive $m$-th root of $1,$ and then by Proposition \ref{Lem2},
\begin{equation*}
\mu _{\Psi ^{m}}(0)>m>\mu .
\end{equation*}%
This contradicts (\ref{k2}). Therefore, $0$ can not be an isolated fixed
point of $\Phi .$
\end{proof}

This corollary is the key ingredient of our method to study invariant
varieties of holomorphic systems. After some analysis about the fixed point
set $Fix(\Phi )$ of $\Phi $ we will be able to prove that $Fix(\Phi )$
contains an analytic variety of complex dimension at least $1,$ consisting
of $0$ and periodic orbits of the same period of the system. Under certain
condition, we will also be able to prove that $Fix(\Phi )$ consists of
periodic orbits of the same period $\frac{2\pi }{|\omega |}.$ Here the
period always means the least positive period.

\section{Periodically \textbf{invariant varieties of holomorphic systems
\label{sec2}}}

\subsection{C\textbf{enters of planar holomorphic systems \label{sec2
copy(1)}}}

In the history of the qualitative theory of planar ODEs, one of the most
interesting problems is to find the condition such that a singularity of an
ODE is a \textit{center} or \textit{isochronous center}.

Let $U$ be a domain in $\mathbb{R}^{2}$ and let $G\in C^{k}(U,\mathbb{R}%
^{2}),k\geq 1.$ For the planar ODE
\begin{equation}
\dot{x}=G(x),x\in U,  \label{I.1}
\end{equation}%
an isolated singular point $p\in U$ is called a \textit{center} if and only
if there exists a punctured neighborhood $V\subset U$ of $p$, consisting of
periodic orbits of (\ref{I.1}) surrounding $p$, and $p$ is called an \textit{%
isochronous center} if and only if it is a center and all orbits of (\ref%
{I.1}) near $p$ have the same period.

Up to now, several classes of systems have been studied extensively, in
relation to the existence of isochronous centers (see [\ref{JS}]). Among
them, the equation
\begin{equation*}
\dot{z}=P(z),z\in U\subset \mathbb{R}^{2}\cong \mathbb{C},
\end{equation*}%
where $P$ is a complex holomorphic function defined on $U$, were considered
in [\ref{BT}], [\ref{CRZ}], [\ref{CD}], [\ref{G1}], [\ref{H1}], [\ref{L}], [%
\ref{Pa}], [\ref{Sa}] and [\ref{V}], etc.

The first result about isochronous centers of complex holomorphic equations
was probably given by Gregor [\ref{G1}] in 1958. Gregor's result can be
concluded as follows. \medskip

\begin{theorem}[Gregor]
\label{th1}\textit{Consider the system}
\begin{equation}
\dot{z}=P(z),z\in U\subset \mathbb{C},  \label{I.2}
\end{equation}%
\textit{where} $P$ \textit{is a holomorphic} \textit{function defined on }$U$%
\textit{. A simple zero }$0\in U$\textit{\ of }$P$\textit{\ is a center of
system (\ref{I.2}) if and only if }$P^{\prime }(0)$\textit{\ is pure
imaginary. In this case }$0$\textit{\ is an isochronous center and the
common period of each cycle surrounding} $0$ \textit{is }$T=\frac{2\pi }{%
|P^{\prime }(0)|}.$
\end{theorem}

This result can be found summarized in [\ref{H1}] and several proofs of this
theorem can be found in [\ref{BT}], [\ref{CRZ}], [\ref{L}], [\ref{Pa}] and [%
\ref{V}].

\subsection{Invariant periodic sets of holomorphic systems in $\mathbb{C}%
^{n} $}

Now, consider an $n$-dimensional complex holomorphic system%
\begin{equation}
\dot{x}=F(x),x\in \Delta ^{n},  \label{I.3}
\end{equation}%
where $\Delta ^{n}$ is a ball centered at the origin $0$ in $\mathbb{C}^{n}$
and $F:\Delta ^{n}\rightarrow \mathbb{C}^{n}$ is a holomorphic mapping such
that $0$ is an isolated zero of $F,$ say, $0$ is an isolated singularity of (%
\ref{I.3}).

When $n=2,$ as a generalization of Gregor's theorem, Needham and McAllister [%
\ref{NM}] proved that if
\begin{equation*}
F^{\prime }(0)=\left(
\begin{tabular}{ll}
$\mu i$ & $0$ \\
$0$ & $\lambda $%
\end{tabular}%
\right) ,
\end{equation*}%
where $\mu \neq 0$ is a real number and $\lambda \neq 0$ is a complex
number, then there exists a one dimensional complex manifold consisting of $%
0 $ and periodic orbits of (\ref{I.3}) of the same period$.$

For arbitrary $n,$ Needham [\ref{N1}] proved the following theorem.

\begin{theorem}[Needham]
If $F^{\prime }(0)$ has a nonzero pure imaginary eigenvalue $\mu i$ and if
the other eigenvalues all have nonzero real parts, then there uniquely
exists a one dimensional complex submanifold of a neighborhood of $0$ in $%
\Delta ^{n}$, consisting of $0$ and periodic orbits of (\ref{I.3}) of the
same period $\frac{2\pi }{|\mu |}$.
\end{theorem}

In this note, we shall use Corollary \ref{C1}, together with some analysis
of fixed point sets of certain moment mapping of the local flow of (\ref{I.3}%
), to prove Needham's theorem in more general versions:

\begin{theorem}
\label{Th1.1}There is a complex analytic variety in $\Delta ^{n}$ with pure
complex dimension at least $1$ consisting of $0$ and periodic orbits of (\ref%
{I.3}) of the same period, if and only if $F^{\prime }(0)$ has a nonzero
pure imaginary eigenvalue.
\end{theorem}

\begin{theorem}
\label{Th1.2}If $F^{\prime }(0)$ has a nonzero pure imaginary eigenvalue $%
\omega i$ and if for any other eigenvalue $\lambda $ of $F^{\prime }(0),\;$%
\begin{equation*}
\lambda /\left( \omega i\right) \neq \pm 2,\pm 3,\dots ,
\end{equation*}%
then there exists a complex analytic variety in $\Delta ^{n}$ of pure
complex dimension at least $1,$ consisting of $0$ and periodic orbits of (%
\ref{I.3}) of the same period $\frac{2\pi }{|\omega |}.$
\end{theorem}

\begin{theorem}
\label{Th1.3}If $F^{\prime }(0)$ has a nonzero pure imaginary eigenvalue $%
\omega i$ and if for any other eigenvalue $\lambda $ of $F^{\prime }(0),$
\begin{equation*}
\lambda /\left( \omega i\right) \neq 0,\pm 1,\pm 2,\dots ,
\end{equation*}%
then there uniquely exists a complex analytic disk $D$ in $\Delta ^{n}$
consisting of $0$ and periodic orbits of (\ref{I.3}) of the same period $%
\frac{2\pi }{|\omega |}.$
\end{theorem}

The term \emph{pure dimension} means that the variety has the same dimension
everywhere, the term \emph{analytic disk} means that $D$ is a holomorphic
embedding of a disk in the complex plane and the uniqueness of $D$ means
that any analytic disk satisfying the condition in the theorem coincides
with $D$ in a neighborhood of $0$ in $\mathbb{C}^{n}.$\medskip

\begin{remark}
The method in this note can also be applied to study the analyticity of
stable manifolds and unstable manifolds of holomorphic differential
equations, by parameter transformation of the time $t$: we can use any line
in the complex plane to replace the real line of $t.$

For example, if we replace the time $t$ by $(a+ib)t$ with $b\neq 0,$ then
the invariant varieties in Theorems \ref{Th1.1} to \ref{Th1.3} become
stable, or unstable, invariant sets in the new systems.
\end{remark}

\section{Proof of Theorems \textbf{\protect\ref{Th1.1}--\protect\ref{Th1.3}}}

\begin{lemma}
\label{lem7}There exist positive numbers $\delta $ and $T_{0},$ such that
for each $T$ in the interval $(0,$ $T_{0}],$ the system (\ref{I.3}) has no
periodic orbit of period $T$ intersecting $\overline{B}_{\delta },$ where $%
B_{\delta }=\{x\in \mathbb{C}^{n};|x|<\delta \}$.
\end{lemma}

\begin{proof}
Let $\phi =\phi (t,x)$ be the local flow of (\ref{I.3})$.$ By Lemma \ref{C2}%
, for any $T_{0}>0,$ there is a $\delta >0,$ such that $\phi (t,x)$ is real
analytic on $[0,T_{0}]\times \overline{B}_{\delta }$ and complex holomorphic
with respect to $x$. We shall show that the conclusion holds for
sufficiently small $\delta $ and sufficiently small $T_{0}$.

Otherwise, for any fixed $\delta $ and $T_{0},$ there exist sequences $%
\{\eta _{j}\}$ in the interval $(0,\delta )$ and $\{t_{j}\}$ in the interval
$(0,T_{0})$ such that
\begin{equation}
\eta _{j}\rightarrow 0\;\mathrm{and\;}t_{j}\rightarrow 0\;\mathrm{as\;}%
j\rightarrow \infty ,  \label{4.1}
\end{equation}%
and that for each $j,$ (\ref{I.3}) has a periodic orbit $\mathcal{O}_{j}$ of
period $t_{j}$ intersecting $B_{\eta _{j}}.$

By the fact that $\phi (0,x)\equiv x$ we have
\begin{equation}
\max_{x\in \overline{B}_{\delta },0\leq t\leq t_{j}}|\phi
(t,x)-x|\rightarrow 0\;\mathrm{as\;}j\rightarrow \infty .  \label{4.2}
\end{equation}

Considering that $\mathcal{O}_{j}\cap B_{\eta _{j}}$ contains infinitely
many fixed points of the time $t_{j}$ mapping $\Phi _{t_{j}}(x)=\phi
(t_{j},x)$ of the local flow $\phi $ and $\eta _{j}<\delta ,$ we conclude by
Lemma \ref{Lem3} that the time $t_{j}$ mapping $\Phi _{t_{j}}$ has a fixed
point $x_{j}\in \partial B_{\delta }$ (the boundary of $B_{\delta }$)$,$ for
each $j.$ Without loss of generality, assume%
\begin{equation}
\lim_{j\rightarrow \infty }x_{j}=x_{0}  \label{4.3}
\end{equation}%
for some $x_{0}\in \partial B_{\delta }.$

We assume that $\delta $ is small enough such that $0$ is the unique
singularity of (\ref{I.3}) in $\overline{B}_{\delta }.$ Then there exist
positive numbers $t_{0}$ and $\theta _{0}$ with $t_{0}<T_{0},$ such that%
\begin{equation*}
|\phi (t_{0},x_{0})-x_{0}|>\theta _{0}.
\end{equation*}%
Thus, by (\ref{4.3}), for sufficiently large $j,$
\begin{equation*}
|\phi (t_{0},x_{j})-x_{j}|>\theta _{0}.
\end{equation*}

On the other hand, since $x_{j}$ is a fixed point of the time $t_{j}$
mapping $\Phi _{t_{j}},$ we have that $\phi (t_{j},x_{j})=x_{j}$ for each $%
j, $ and then, putting $t_{0}=m_{j}t_{j}+r_{j}$ for each $j,$ where $m_{j}$
is an integer and $0\leq r_{j}<t_{j},$ we have%
\begin{eqnarray*}
&&|\phi (r_{j},x_{j})-x_{j}|=|\phi (r_{j},\phi (m_{j}t_{j},x_{j}))-x_{j}| \\
&=&|\phi (m_{j}t_{j}+r_{j},x_{j})-x_{j}|=|\phi (t_{0},x_{j})-x_{j}|>\theta
_{0}.
\end{eqnarray*}%
This contradicts (\ref{4.2}), for $0\leq r_{j}<t_{j}$ and $x_{j}\in \partial
B_{\delta }.$ This completes the proof.
\end{proof}

To prove Theorems \ref{Th1.1}--\ref{Th1.3}, we need some known results about
analytic sets. The reader is referred to [\ref{Ch}], pages 14--57, for the
details.\label{ddd}

\begin{lemma}
\label{cr1}Let $\mathcal{A}$ be an analytic subset of an open subset of $%
\mathbb{C}^{n}.$ Then for any $p\in \mathcal{A},$ the (complex) dimension $%
\dim _{p}\mathcal{A}$ of $\mathcal{A}$ at $p$ is at least $1$ if and only if
$p$ is an accumulated point of $\mathcal{A}.$
\end{lemma}

\begin{lemma}
\label{cr2}If $\dim _{p}\mathcal{A}\geq 1,$ then $\mathcal{A}$ has an
irreducible component $\mathcal{B}$ containing $p,$ with a pure complex
dimension at least $1.$
\end{lemma}

\begin{lemma}
\label{cr4}Any irreducible component of $\mathcal{A}$ is the closure in $%
\mathcal{A}$ of some connected component of $\mathrm{reg}\mathcal{A}$, where
$\mathrm{reg}\mathcal{A}$ is the set of all regular points of $\mathcal{A}.$
\end{lemma}

A point $p\in \mathcal{A}$ is called regular, if there is a neighborhood $%
V_{p}$ of $p\ $in $\mathbb{C}^{n}$, such that $V_{p}\cap \mathcal{A}$ is a
complex submanifold of $V_{p}.$

\begin{proof}[\textbf{PROOF OF THEOREM} \textbf{\protect\ref{Th1.1}}.]
We first proof the necessity. Assume that there exists an analytic variety $%
\Sigma $ of pure complex dimension at least $1,$ consisting of the
singularity $0$ and periodic orbits of (\ref{I.3}) of the same period $\tau
_{0}$ (note that the period always indicates the least positive period$.$)

By Lemma \ref{lem7}, there exist positive numbers $\delta $ and $T_{0}$ $%
(<\tau _{0})$ such that for each $T\in (0,$ $T_{0}],$ the system (\ref{I.3})
has no periodic orbit of period $T$ intersecting $\overline{B}_{\delta }$.

Let $M$ be a prime number such that $\tau _{0}/M<T_{0}$ and let $\Phi
(x)=\phi (\tau _{0}/M,x)$ be the time $\tau _{0}/M\;$mapping of the local
flow $\phi (t,x)$ of (\ref{I.3})$.$ Then $0$ is the unique fixed point of $%
\Phi $ located in $\overline{B}_{\delta }.$ On the other hand, by the
assumption about $\Sigma ,$ $0$ is an accumulated point of $\Sigma $ and all
points of $\Sigma $ are fixed points of $\Phi ^{M}(x)=\phi (\tau _{0},x),$
and then $0$ is an accumulated point of fixed points of $\Phi ^{M}.$

Thus by Lemma \ref{lem6}, $\Phi ^{\prime }(0)$ has an eigenvalue $\Lambda $
such that $\Lambda \neq 1$ but $\Lambda ^{m}=1,$ which implies by Lemma \ref%
{C2} that $F^{\prime }(0)$ has an eigenvalue $\lambda $ with $\Lambda =e^{%
\frac{\tau _{0}\lambda }{M}}\neq 1$ and $\Lambda ^{M}=e^{\tau _{0}\lambda
}=1.$ Clearly, $\lambda =\omega i$ fore some real number $\omega \neq 0$.
This completes the proof of the necessity.

The sufficiency follows from Theorem \ref{Th1.2} directly.
\end{proof}

\begin{proof}[\textbf{PROOF OF\quad THEOREM \protect\ref{Th1.2}.}]
Assume $F^{\prime }(0)$ has an eigenvalue $\omega i$ such that $\omega \neq
0 $ is a real number and for any other eigenvalue $\lambda $%
\begin{equation}
\lambda /\left( \omega i\right) \neq \pm 2,\pm 3,...  \label{cond}
\end{equation}%
Then by Corollary \ref{C1}, the origin $0$\ is an accumulated point of fixed
points of the time $\frac{2\pi }{|\omega |}$\ mapping $\Phi =\phi (\frac{%
2\pi }{|\omega |},\cdot )$\ of the local flow $\phi (t,x)$\ of (\ref{I.3})\
defined in a neighborhood of the origin.

By Lemma \ref{C2}, for some ball $B$ centered at the origin $0,$ $\phi (t,x)$
is real analytic on $[0,\frac{2\pi }{|\omega |}]\times \overline{B}$ and
holomorphic with respect to $x\in \overline{B}.$ Let $\mathcal{A}_{B}$ be
the set of all fixed points of the mapping $\Phi $ in
\begin{equation*}
V_{B}=\phi ([0,\frac{2\pi }{|\omega |}]\times B)=\{\phi (t,x);t\in \lbrack 0,%
\frac{2\pi }{|\omega |}],x\in B\}.
\end{equation*}%
Then $\mathcal{A}_{B}$ is an analytic subset of $V_{B}$ and by Lemma \ref%
{cr1} we have $\dim _{0}\mathcal{A}_{B}\geq 1$, and then, by Lemma \ref{cr2}%
, there exists an irreducible component $\Sigma _{B}$ of $\mathcal{A}_{B}$
containing $p$ with a pure complex dimension at least $1$.

We first show that $\Sigma _{B}\backslash \{0\}$ consists of periodic orbits
of (\ref{I.3}). It suffices to prove that $\Sigma _{B}$ is invariant by the
flow $\phi (t,x).$ Note that by the definition of $V_{B},$ $\mathcal{A}_{B}$
is invariant by the local flow $\phi ,$ say, $\phi (t,\mathcal{A}_{B})=%
\mathcal{A}_{B}$ for all real $t.$

By Lemma \ref{cr4}, there exists a connected component $S$ of \textrm{Reg}$%
\mathcal{A}_{B}$ such that $\Sigma _{B}$ is the closure of $S$ in $\mathcal{A%
}_{B}$. Then each $p\in S$ has a neighborhood $U_{p}$ in $V_{B}$, such that $%
S\cap U_{p}$ is a complex submanifold of $U_{p}$ and $S$ is the only
connected component of $\mathcal{A}_{B}$ intersecting $U_{p}$, and then $%
\Sigma _{B}$ is the only irreducible component of $\mathcal{A}_{B}$
intersecting $U_{p}$. On the other hand, for sufficiently small $\varepsilon
>0,\Phi _{\varepsilon }(\Sigma _{B})=\phi (\varepsilon ,\Sigma _{B})$ must
intersects $U_{p},$ and then by the fact that $\Phi _{\varepsilon }=\phi
(\varepsilon ,\cdot )$ is a biholomorphic mapping from $V_{B}$ onto $\Phi
_{\varepsilon }(V_{B}),$ $\Phi _{\varepsilon }(\Sigma _{B})$ is also an
irreducible component of $\Phi _{\varepsilon }(\mathcal{A}_{B})=\mathcal{A}%
_{B}$ (Recall that $\mathcal{A}_{B}$ is invariant by the flow $\phi (t,x)).$
Therefore, $\Phi _{\varepsilon }(\Sigma _{B})$ and $\Sigma _{B}$ are both
irreducible components of $\mathcal{A}_{B}$ and both intersect $U_{p},$ and
then $\Phi _{\varepsilon }(\Sigma _{B})=\Sigma _{B},$ which implies that $%
\phi (t,\Sigma _{B})=\Sigma _{B}$ for all real $t.$ Thus $\Sigma _{B}$ is
invariant by $\phi $ and then $\Sigma _{B}\backslash \{0\}$ consists of
periodic orbits.

It is clear that for each periodic orbit $\mathcal{O}\subset \Sigma
_{B}\backslash \{0\}$, the period $T(\mathcal{O})$ of $\mathcal{O}$ equals
to $\frac{2\pi }{|\omega |m_{\mathcal{O}}}$ for some positive integer $m_{%
\mathcal{O}}$. By Lemma \ref{lem7}, we may assume that $B$ is small enough
such that for some fixed $T_{0}>0$, the period of each orbit in $\Sigma
_{B}\backslash \{0\}$ is larger than $T_{0}$. Then for each periodic orbit $%
\mathcal{O}\subset \Sigma _{B}$ we have $T_{0}\leq T(\mathcal{O})=\frac{2\pi
}{|\omega |m_{\mathcal{O}}},$ with
\begin{equation}
m_{\mathcal{O}}\in \{1,2,\dots ,m^{\ast }\},  \label{4.5}
\end{equation}%
where $m^{\ast }$ is the integral part of $2\pi \left( T_{0}|\omega |\right)
^{-1}$.

We assert that if the ball $B$ is small enough, then for each periodic orbit
$\mathcal{O}\subset \Sigma _{B}\backslash \{0\},$ $m_{\mathcal{O}}=1.$
Otherwise, by (\ref{4.5}) there is a fixed integer $m>1,$ and a sequence $%
\{x_{k}\}\subset \Sigma _{B}\backslash \{0\},$ such that $x_{k}\rightarrow 0$
as $k\rightarrow \infty $ and the periodic orbit passing through $x_{k}$ has
period $\frac{2\pi }{|\omega |m}$ for each $k,$ say, $\phi (\frac{2\pi }{%
|\omega |m},x_{k})=x_{k}$ for each $k.$

Let $M$ be any given prime number with $\frac{2\pi }{|\omega |mM}<T_{0}.$
Then $0$ is an isolated fixed point of $\Theta (x)=\phi (\frac{2\pi }{%
|\omega |mM},x)$ but $0$ is an accumulated point of fixed points of $\Theta
^{M},$ for
\begin{equation*}
\Theta ^{M}(x_{k})=\phi (\frac{2\pi }{|\omega |m},x_{k})=x_{k},k=1,2,\dots
\end{equation*}%
Therefore, by Lemma \ref{lem6}, the Jacobian matrix $\Theta ^{\prime }(0)$
must have an eigenvalue $\Lambda $ with $\Lambda \neq 1$ but $\Lambda
^{M}=1. $ Hence, by (\ref{k1}), we have $\Lambda =e^{\frac{2\pi }{|\omega |mM%
}\lambda }$ for some eigenvalue $\lambda $ of $F^{\prime }(0)$ such that $%
\lambda \neq 0$ and $\frac{2\pi }{|\omega |m}\lambda $ is a multiple of $%
2\pi i,$ and then $\lambda /(\omega i)=\pm km$ for some positive integer $k,$
which contradicts (\ref{cond}). Therefore, we have $m=1.$ The proof is
complete.
\end{proof}

\begin{proof}[\textbf{PROOF OF\ THEOREM }\protect\ref{Th1.3}.]
Assume that $\lambda _{1},\dots ,\lambda _{n}$ are the $n$ eigenvalues of
the Jacobian matrix $F^{\prime }(0)$ of $F$ at $0,$ with $\lambda
_{1}=\omega i$ and
\begin{equation}
\lambda _{l}/\lambda _{1}\neq 0,\pm 1,\pm 2,\pm 3,...,l=2,3,\dots ,n.
\label{a2}
\end{equation}%
Ignoring a linear transform of the phase space, we may assume that the
Jacobian matrix $F^{\prime }(0)$ is lower triangular and has $\lambda
_{1},\lambda _{2},...,\lambda _{n}$ down its main diagonal. We write this by%
\begin{equation}
F^{\prime }(0)=\left( \lambda _{1},\lambda _{2},\lambda _{3},...,\lambda
_{n}\right) .  \label{a1}
\end{equation}

By Lemma \ref{C2}, there exists a ball $B$ centered at the origin $0$ such
that the local flow $\phi (t,x)$ of (\ref{I.3}) is real analytic on $[0,%
\frac{2\pi }{|\omega |}]\times \overline{B}$ and complex holomorphic with
respect to $x\in B.$

By (\ref{a1}) and Lemma \ref{C2},
\begin{equation*}
\Lambda _{l}=e^{\frac{2\pi }{|\omega |}\lambda _{l}},l=1,2,\dots ,n,\;
\end{equation*}%
are all the eigenvalues of the Jacobian matrix $\Phi ^{\prime }(0)$ of the
time $\frac{2\pi }{|\omega |}$ mapping $\Phi (x)=\phi (\frac{2\pi }{|\omega |%
},x)$ at $0,$ and by (\ref{a2}) we have
\begin{equation}
\Lambda _{l}\neq 1,l=2,\dots ,n.  \label{a3}
\end{equation}

It also follows from (\ref{a1}) and Lemma \ref{C2} that $\Phi ^{\prime
}(0)=\left( 1,\Lambda _{2},\Lambda _{3}\dots ,\Lambda _{n}\right) $ is a
lower triangular matrix which has $1,\Lambda _{2},\Lambda _{3}\dots ,\Lambda
_{n}$ down its main diagonal. Therefore, putting $x=(x_{1},x_{2},\dots
,x_{n})$ and $\Phi =(\varphi _{1},\varphi _{2},\dots ,\varphi _{n}),$ by (%
\ref{a3}) we have
\begin{equation*}
\left. \det \left( \frac{\partial \left( \varphi _{2},\varphi _{3},\dots
,\varphi _{n}\right) }{\partial \left( x_{2},x_{3},\dots ,x_{n}\right) }%
-I_{n-1}\right) \right\vert _{x=0}\neq 0
\end{equation*}%
where $I_{n-1}$ is the $(n-1)\times (n-1)$ unit matrix. So, by the implicit
function theorem, there uniquely exist one variable complex holomorphic
functions $x_{l}=x_{l}(x_{1}),l=2,\dots ,n,$ defined in a neighborhood of
the origin in the complex plane $\mathbb{C},$ solving the system of the
equations
\begin{equation*}
\varphi _{l}(x_{1},x_{2},\dots ,x_{n})=x_{l},l=2,\dots ,n,
\end{equation*}%
in a neighborhood of the origin $0,$ with
\begin{equation}
x_{l}(0)=0,l=2,3,\dots ,n.  \label{a7}
\end{equation}%
Therefore, in a neighborhood of the origin $0$ in $\mathbb{C}^{n},$ the
fixed point equation
\begin{equation}
\Phi (x_{1},x_{2},\dots ,x_{n})=(x_{1},x_{2},\dots ,x_{n})  \label{a4}
\end{equation}%
is equivalent to the system of the $n$ equations
\begin{equation}
\left\{
\begin{tabular}{l}
$x_{1}=\varphi _{1}(x_{1},x_{2}(x_{1}),\dots ,x_{n}(x_{1})),$ \\
$x_{l}=x_{l}(x_{1}),l=2,\dots ,n.$%
\end{tabular}%
\ \right.  \label{a5}
\end{equation}

By Corollary \ref{C1}, $0$ is an accumulated point of fixed points of $\Phi
, $ and then $0$ is an accumulated point of zeros of the holomorphic
function $x_{1}-\varphi _{1}(x_{1},x_{2}(x_{1}),\dots ,x_{n}(x_{1}))$ which
is defined in a neighborhood of the origin $0$ in $\mathbb{C}.$ Therefore,
we have
\begin{equation*}
\varphi _{1}(x_{1},x_{2}(x_{1}),\dots ,x_{n}(x_{1}))\equiv x_{1}
\end{equation*}%
in a neighborhood of $x_{1}=0$ in $\mathbb{C}$, and then the equation (\ref%
{a5}) reads%
\begin{equation}
\left\{
\begin{tabular}{l}
$x_{1}=x_{1},$ \\
$x_{l}=x_{l}(x_{1}),l=2,\dots ,n.$%
\end{tabular}%
\right.  \label{a6}
\end{equation}%
It is clear that there exists a positive number $\delta $ such that
\begin{equation*}
\Sigma _{\delta }^{\ast }=\{(x_{1},\varphi _{2}(x_{1}),\dots ,\varphi
_{n}(x_{1}));\;x_{1}\in \mathbb{C},|x_{1}|<\delta \}
\end{equation*}%
is an analytic disk in $\Delta ^{n},$ which contains $0$ by (\ref{a7}).

By the equivalence of (\ref{a4}) and (\ref{a6}), in a neighborhood of $0,$ $%
\Sigma _{\delta }^{\ast }$ coincides with the set $\mathcal{A}_{B}$ of all
fixed points of $\Phi $ in
\begin{equation*}
V_{B}=\phi ([0,\frac{2\pi }{|\omega |}]\times B),
\end{equation*}%
and then $0\in \mathcal{A}_{B}$ has a neighborhood $U_{0}$ in $\mathcal{A}%
_{B}$ so that $U_{0}$ is an analytic disk in $\mathbb{C}^{n}.$ On the other
hand, by Theorem \ref{Th1.2} there exists an analytic variety $\Sigma $ of
pure complex dimension at least $1,$ consisting of $0$ and periodic orbits
of the same period $\frac{2\pi }{|\omega |}$ of the system (\ref{I.3}), and
then by the equivalence of (\ref{a4}) and (\ref{a6}), $0\in \Sigma $ has a
neighborhood $V_{0}$ in $\Sigma $ such that $V_{0}\subset U_{0},$ which
implies that the dimension of $\Sigma $ is equal to $1,$ and then $V_{0}$
coincides with $U_{0}\mathcal{\ }$in a neighborhood of $0.$

Now, we can conclude that $\Sigma ,\Sigma _{\delta }^{\ast }$ and $\mathcal{A%
}_{B}$ coincide in a neighborhood of the origin. Therefore, when $B$ is
small enough, $\mathcal{A}_{B}$ is contained in $\Sigma \cap \Sigma _{\delta
}^{\ast },$ and then $\mathcal{A}_{B}$ contains an analytic disk $D$
consisting of $0$ and periodic orbits of the same period $\frac{2\pi }{%
|\omega |}$.

The above argument also show that any analytic disk $D$ consisting of $0$
and periodic orbits of (\ref{I.3}) with the same period $\frac{2\pi }{%
|\omega |}$ coincides with $\Sigma _{\delta }^{\ast }$ in a neighborhood of $%
0,$ which implies the uniqueness and the proof is complete.
\end{proof}

\section{Appendix}

\begin{proof}[Proof of Proposition \protect\ref{Lem2}]
By the given condition, there exists a ball $B\subset \Delta ^{n}$ with
center $0$ such that both $f$ and $f^{m}$ is well defined in a neighborhood
of $\overline{B}$ and $0$ is the unique fixed point of both $f$ and $f^{m}$
in $\overline{B}.$

Let $\lambda $ be an eigenvalue of $f^{\prime }(0)$ that is a primitive $m$%
-th root of $1.$ Then the Jacobian matrix $(f^{m})^{\prime }(0)=(f^{\prime
}(0))^{m}$ of $f^{m}$ at $0$ has the eigenvalue $\lambda ^{m}=1,$ and then
by Lemma \ref{lem4} we have%
\begin{equation}
\mu _{f^{m}}(0)\geq 2.  \label{6}
\end{equation}

We first assume that none of eigenvalues of $f^{\prime }(0)$ equals to $1.$
Then applying Lemma \ref{lem4} once more we have
\begin{equation}
\mu _{f}(0)=1.  \label{5}
\end{equation}

We can construct a sequence $f_{k}:\Delta ^{n}\rightarrow \mathbb{C}^{n}$ of
holomorphic mappings$,$ converging to $f$ uniformly, such that for each $k,$
$0$ is a simple fixed point of both $f_{k}$ and $f_{k}^{m},$ say, $%
f_{k}(0)=f_{k}^{m}(0)=0,$ but neither $f_{k}^{\prime }(0)$ nor $%
(f_{k}^{m})^{\prime }(0)$ has eigenvalue $1.$ This can be accomplished by
small perturbations of the linear part of $f$ at $0$.

It is clear that for sufficiently large $k,$ $f_{k}^{m}$ is well defined on $%
\overline{B}$ and $f_{k}^{m}$ converges to $f^{m}$ uniformly on $\overline{B}
$ as $k$ tends to $\infty .$ Therefore, by (\ref{6}), Lemma \ref{lem5} and
the condition that $0$ is a simple fixed point of $f_{k}^{m}$, we conclude
that for sufficiently large $k,$ $f_{k}^{m}$ has another fixed point $%
x_{k}\neq 0$ in $B.$ On the other hand, since $0$ is the unique fixed point
of $f$ in $\overline{B}$, by Lemma \ref{lem5} and (\ref{5}), for
sufficiently large $k,$ $0$ is the unique fixed point of $f_{k}$ in $B.$ Thus%
$,$ we have $f_{k}(x_{k})\neq x_{k}.$

Now that $f_{k}(x_{k})\neq x_{k},f_{k}^{m}(x_{k})=x_{k}$ and $m$ is prime, $%
x_{k}$ is a periodic point of $f_{k}$ with period $m,$ say, the periodic
orbit $\Gamma (x_{k})=\{x_{k},f_{k}(x_{k}),\dots ,f_{k}^{m-1}(x_{k})\}$
contains $m$ distinct points. Since $0$ is the unique fixed point of $f^{m}$
in $\overline{B}$, by the convergence of $f_{k}^{m}$ we have that $%
x_{k}\rightarrow 0$ as $k\rightarrow \infty ,$ and then by the condition $%
f(0)=0$ and the convergence of $f_{k},$ the periodic orbit $\Gamma (x_{k})$
is contained in $B$ for sufficiently large $k.$ Thus for sufficiently large $%
k,$ $0,x_{k},f_{k}(x_{k}),\dots ,f_{k}^{m-1}(x_{k})$ are $m+1$ distinct
fixed points of $f_{k}$ in $B.$ By Lemma \ref{lem5}, we then have%
\begin{equation*}
\mu _{f^{m}}(0)\geq m+1.
\end{equation*}%
Therefore, (\ref{3}) holds.

In general, we can consider the mapping%
\begin{equation*}
f_{\varepsilon }=(\varepsilon _{1}z_{1},\dots ,\varepsilon
_{n}z_{n})+f(z_{1},\dots ,z_{n})
\end{equation*}%
where $\varepsilon =(\varepsilon _{1},\dots ,\varepsilon _{n})$ is so chosen
that $f_{\varepsilon }^{\prime }(0)$ has eigenvalue $\lambda $ but none
other eigenvalue equals to $1.$ It is easy to see that we can chose $%
\varepsilon $ sufficiently small so that $f_{\varepsilon }^{m}$ is
sufficiently close to $f^{m},$ and then by Lemma \ref{lem5}, we have $\mu
_{f^{m}}(0)\geq \mu _{f_{\varepsilon }^{m}}(0)\geq m+1.$ This completes the
proof.
\end{proof}

\end{document}